\documentclass[reqno,11pt]{amsart}

\newtheorem{theorem}{Theorem}[section]
\newtheorem{lemma}[theorem]{Lemma}

\theoremstyle{definition}

\newtheorem{assumption}{Assumption}[section]

 \theoremstyle{remark}
\newtheorem{remark}[theorem]{Remark}

\newcommand\bH{\mathbb{H}}
\newcommand\bL{\mathbb{L}}

\newcommand\bR{\mathbb{R}}

\newcommand\cF{\mathcal{F}}
\newcommand\cH{\mathcal{H}}
 
\newcommand\cP{\mathcal{P}}

\newcommand{\nlimsup}{\operatornamewithlimits{\overline{lim}\,}}

\newcommand{\mysection}[1]{\section{#1}
\setcounter{equation}{0}}

\newcommand\cbrk{\text{$]$\kern-.15em$]$}} 
\newcommand\opar{
\text{\,\raise.2ex\hbox{${\scriptstyle |}$}\kern-.34em$($}} 
\newcommand\cpar{%
\text{$)$\kern-.34em\raise.2ex\hbox{${\scriptstyle |}$}}\,}
\newcommand\obrk{\text{$[$\kern-.15em$[$}}

\begin{document}

\title[Filtering equations with Lipschitz
  coefficients]
{Filtering equations for partially
observable diffusion processes with Lipschitz
continuous coefficients}
\author{N.V. Krylov}
\address{127 Vincent Hall, University of Minnesota, Minneapolis,
 MN, 55455}
\thanks{The work   was partially supported by
NSF Grant DMS-0653121}
\email{krylov@math.umn.edu}

\renewcommand{\subjclassname}{%
\textup{2000} Mathematics Subject Classification}


\begin{abstract}
We present several results on   smoothness in $L_{p}$ sense
of filtering densities under the Lipschitz continuity  assumption  
on the coefficients of a partially observable diffusion processes.
We obtain them by rewriting in divergence form
  filtering equation
which are usually considered in terms of formally adjoint to
operators  in nondivergence form.
 
\end{abstract}

\maketitle
\mysection{Introduction}
                                     \label{section 4.17.1}

For the author, one of the main motivations for
developing the theory of stochastic partial differential
equations (SPDEs) is its relation to the filtering
problem for partially observable diffusion processes.

This problem's setting is as follows.

Let $(\Omega,\cF,P)$ be a complete probability space
with an increasing filtration $\{\cF_{t},t\geq0\}$
of complete, with respect to $(\cF,P)$, $\sigma$-fields
$\cF_{t}\subset\cF$. Denote by $\cP$ the predictable
$\sigma$-field in $\Omega\times(0,\infty)$
associated with $\{\cF_{t}\}$. Let $d\geq1$, $d_{1}>d$, 
and $d_{2}\geq d_{1}$
be integers and $w_{t}$ be a $d_{2}$-dimensional
Wiener process with respect to $\{\cF_{t}\}$. Let $K,T,
\delta>0$ be fixed finite constants.

Consider a $d_{1}$-dimensional two component process
 $z_{t}=(x_{t},y_{t})$
with $x_{t}$ being $d$-dimensional and $y_{t}$ 
$(d_{1}-d)$-dimensional. We
assume that $z_{t}$ is a diffusion process
 defined as a solution of the system
\begin{equation}\begin{split}         
                                            \label{eq3.2.14} 
& dx_{t}=b(t,z_{t}) dt+\theta (t,z_{t})dw_{t}, \\ 
& dy_{t}=B(t,z_{t}) dt+\Theta (t,y_{t})dw_{t}
\end{split}\end{equation}
with some initial data.

The coefficients of \eqref{eq3.2.14} are assumed to be
 vector- or matrix-valued
functions of appropriate dimensions defined on
 $[0,T]\times\bR^{d_{1}}$.
Actually $\Theta(t,y)$ is assumed to be independent
 of $x$, so  that it
is a function on $[0,T]\times\bR^{d_{1}-d}$ rather than 
$[0,T]\times\bR^{d_{1}}$
but as always we may think of $\Theta(t,y)$ as a
 function of $(t,z)$ as well.

The component $x_{t}$ is treated as unobservable
and $y_{t}$ as the only observations available.
The problem is to find a way to compute the density
$\pi_{t}(x)$
of the conditional 
distribution of $x_{t}$ given $y_{s},s\leq t$.
Finding an equation satisfied by $\pi_{t}$
(filtering equation) is considered
to be a solution of the (filtering) problem.
Filtering equations turn out to be 
particular cases of SPDEs.

The history of filtering equations for diffusion 
processes is long and its beginning is
controversial.
Probably, the first filtering equations 
were published in \cite{St}. They turned out to be plain wrong.
Then in \cite{Ku64} other equations were proposed, see for instance
equation (5) of \cite{Ku64}. However, it is hard to make
sense of these equations because most likely some terms
appeared from stochastic integrals written in
 the Stratonovich form and the others appeared from
  the It\^o integrals. Perhaps, the author
of \cite{Ku64} realized this too and published
an attempt to rescue some results of \cite{Ku64}
 in \cite{Ku67}. This attempt
  turned  successful for  simplified models
  without the so-called cross terms.

Meanwhile, in \cite{Sh66} the correct filtering
 equations in full generality,
yet assuming some regularity of the filtering density,
 were presented and then in \cite{LS68}
they were rigorously proved.
This is the reason we propose to call the filtering equations
in the case of partially observable diffusion
processes {\em Shiryaev's equations\/}
and their particular case without cross terms
{\em Kushner's equations\/}.

In case $d=1$ the result of \cite{Sh66} is presented
 in \cite{LS}
on the basis of  
the famous Fujisaki-Kallianpur-Kunita theorem (see \cite{FKK})
about the filtering equations in a very general setting
(much more general than in   \cite{LS68}). Some authors
even call the filtering equation for diffusion processes
the Fujisaki-Kallianpur-Kunita equation.

By adding to the Fujisaki-Kallianpur-Kunita theorem
 some simple facts from the theory of SPDEs,
the a priori regularity assumption  was removed 
  in \cite{KR} and under the Lipschitz  
and uniform nondegeneracy assumption
the $L_{2}$-version
of Theorem \ref{thm3.2.22} was proved. 
The basic result of \cite{KR} is that $\pi_{t}\in H^{1}_{2}$.
It is also proved that
if the coefficients are smoother, $\pi_{t}(x)$ is smoother too.
The nondegeneracy assumption is removed 
in \cite{R} on the account of assuming
 that $\theta\theta^{*}$ is three
times continuously differentiable in $x$. It is again
proved that $\pi_{t}\in H^{1}_{2}$ and $\pi_{t}$ is even smoother
if the coefficients are smoother.

In \cite{K99} the results of \cite{KR} were improved,
$\theta\theta^{*}$ is assumed to be twice
  continuously differentiable in $x$ and it is shown that
$\pi_{t}\in H^{2}_{p}$ with any $p\geq2$.
 
The above mentioned results of \cite{KR}, \cite{R},
and \cite{K99}
use the filtering theory
in combination with the theory of SPDEs, the latter
being  stimulated by certain needs of  filtering theory. 
It turns out  that the theory of SPDEs alone can be used
to obtain the above mentioned regularity
 results about $\pi_{t}$ without knowing anything from 
the filtering theory itself. It also can be used to solve
other problems from the filtering theory.

The first ``direct'' (only using the theory of SPDEs)
proof of regularity of $\pi_{t}$
is given in \cite{KZ} in the case that
  system
\eqref{eq3.2.14} defines a nondegenerate diffusion process
and $\theta\theta^{*}$ is twice
 continuously differentiable in $x$. 
It is proved that $\pi_{t}\in H^{2}_{p}$ with any $p\geq2$
as in \cite{K99}.
Advantages of having arbitrary $p$ are seen from results
like our Theorem \ref{theorem 1.19.1}.
Of course, on the way of investigating $\pi_{t}$
in \cite{KZ} filtering equations are derived 
``directly'' in an absolutely
different manner than before (on the basis of an idea
from \cite{KR81}).

In this article we relax the smoothness assumption 
in \cite{KZ} to the assumption that the coefficients
of \eqref{eq3.2.14} are merely Lipschitz continuous,
the assumption which is almost always supposed to hold
when one deals with systems like \eqref{eq3.2.14}.
We find that $\pi_{t}\in H^{1}_{p}$.
Thus, under the weakest smoothness assumptions we obtain
the best (in the author's opinion)
 regularity result on $\pi_{t}$. 
In particular, we prove that if the initial data
is sufficiently regular, then the filtering density
is almost Lipschitz continuous in $x$ and $1/2$
H\"older continuous in $t$.
However, we still
assume $z_{t}$ to be nondegenerate.   
Our approach is heavily based on analytic results.
There is also a probabilistic approach developed in \cite{Ku2}
and based on explicit formulas for solutions initiated
in \cite{Pa} and later developed
in \cite{KR81}  and \cite{Ku1} (also see references therein). This approach cannot give
as sharp results as ours in our situation.  

It seems to the author that 
under the same assumptions of Lipschitz continuity,
 by following an idea
from 
  \cite{K79} one can solve another problem from filtering theory,
  the so-called innovation problem, and obtain the equality
$$\sigma\{y_{s},s\leq t\}
=\sigma\{\check{w}_{s},s\leq t\},
$$
where $\check{w}_{t}$ is the innovation Wiener process
of the problem (its definition is reminded
in Section \ref{section 2.2.1}). Recall that
for degenerate diffusion processes the positive solution
of the innovation problem is obtained in \cite{Pu}
again on the basis of the theory of SPDEs under the assumption
that the coefficients are more regular.

By the way, in our situation, if the coefficients are more regular,
the filtering
equation can be rewritten in a nondivergence form and then
additional smoothness of the filtering density,
existence of which is already established in this article,
is obtained on the basis of regularity results from \cite{K99}.

The article is organized as follows. In Section \ref{section 2.2.1}
we state our main results part of which is 
proved in the same section. In Sections \ref{section 1.31.1}
and \ref{section 1.31.2}
we prove Theorems \ref{thm3.2.22}  
 and \ref{theorem 1.28.1}, respectively. 
Section \ref{section 1.31.3} contains
a collection of results from the theory of SPDEs
which we use in  
the previous sections.

As it is done traditionally in filtering theory
we consider finite-dimensional driving Wiener
processes. However, our results will be based
on the theory of SPDEs,
outlined in Section \ref{section 1.31.3}, with countably many Wiener processes.
We leave to the reader to do
some trivial modifications in Section \ref{section 1.31.3}
 in order to be able
to apply its results in such cases.

\mysection{Main results}
                                         \label{section 2.2.1}
First we state and discuss our assumptions.

\begin{assumption}
                                             \label{asm3.2.15}
The functions $b$, $\theta $, $B$, and $\Theta $ are
 Borel measurable and bounded
functions of their arguments. Each of them satisfies 
the Lipschitz
condition in $z$ with constant $K\in(0,\infty)$. 
\end{assumption}

Introduce
\begin{equation}                       \label{eq3.2.19.4}
\tilde{\theta }_{t}(z) =\left(
\begin{array}{cc}
\theta(t,z)\\
\Theta(t,y)
\end{array}
\right),\quad
\tilde{a}_{t}(z) =\frac{1}{2}
\tilde{\theta }_{ t}\tilde{\theta }^{*}_{t}(z),\quad
\tilde{b}_{t}(z) =\left(
\begin{array}{cc}
b(t,z) \\
B(t,z)
\end{array}\right),
\end{equation}
\begin{equation}\label{eq3.2.19.1}
\tilde{L}_{t}(z)  =   \tilde{a}^{ij}_{t}(z)
\frac{\partial^{2} }{\partial z^{i} \partial z^{j}}+
 \tilde{b}^{i}_{t}(z)\frac{\partial }{\partial z^{i}},
\end{equation}
where $\tilde{\theta}^{*}$ is the transpose of $\tilde{\theta}$
and the summation convention is imposed.

\begin{remark}                         \label{rem3.2.19.2.0}
System of equations \eqref{eq3.2.14} can be now written as
\begin{equation}
\label{10.13.2}
dz_{t}=\tilde{b}(t,z_{t})dt+\tilde{\theta}(t,z_{t})dw_{t}.
\end{equation}
\end{remark}

\begin{assumption}
                                      \label{assumption 1.15.1} 
The process $z_{t}$ is uniformly nondegenerate:
for any $\lambda,z\in\bR^{d_{1}}$ and $t\in[0,T]$ we have
$$
 \tilde{a}^{ij}_{t}(z)\lambda^{i}
\lambda^{j}\geq\delta|\lambda|^{2}.
$$
\end{assumption}

Traditionally, Assumption \ref{assumption 1.15.1}
is split into two following assumptions in which some
useful objects are introduced. These assumptions
 were also used in the past to reduce $\tilde{\theta}$
to the so-called triangular form by replacing
$w_{t}$ with a different Brownian motion.

\begin{assumption}
                                           \label{asm3.2.16}
The symmetric matrix $\Theta \Theta^{*}$ is 
invertible and 
$$
\Psi :=(\Theta
\Theta^{*} )^{-\frac{1}{2}}
$$ 
is a bounded function of $(t,y)$.
\end{assumption}
\begin{remark}
                                           \label{rem3.2.16.1}
Assumption~\ref{asm3.2.16} follows
from Assumption \ref{assumption 1.15.1}
and, furthermore, $\Psi\leq\delta^{-1}(\delta^{ij})$.
\end{remark}
\begin{assumption}\label{asm3.2.17}
For any $\xi \in \bR^{d}$, $z=(x,y)\in \bR^{d_{1}}$, and $t>0$,
we have
\[
|Q(t,y)\theta^{*}(t,z)\xi |^{2}\geq \delta |\xi |^{2},
\]
where $Q$ is the orthogonal projector on $\text{Ker}\,\Theta$. 
In other words,
\begin{equation}                        \label{eq3.2.17.1}
(\theta (I-\Theta^{*}\Psi^{2}\Theta )
\theta^{*} \xi ,\xi )\geq 
\delta |\xi|^{2}.
\end{equation}
\end{assumption}
\begin{remark}\label{rem3.2.17.2}
From \eqref{eq3.2.17.1} we see that $\theta \theta^{*}$ is
uniformly positive definite with constant of positivity $\delta$.
Also, it turns out that
\eqref{eq3.2.17.1}
  holds under Assumption~\ref{assumption 1.15.1}.

Indeed, take a $\zeta=(\xi,\Psi\eta)\in
\bR^{d}\times \bR^{d_{1}-d}$ with
$\eta=-\Psi\Theta\theta^{*}\xi$ and observe that
$$
2\delta|\xi|^{2}\leq2(\tilde{a}\zeta,\zeta)=
|\tilde{\theta}^{*}\zeta|^{2}
=|\theta\xi|^{2}+2(\tilde{\theta}^{*}\xi,\Theta^{*}\Psi\eta)
+|\Theta^{*}\Psi\eta|^{2}
$$
$$
=|\theta\xi|^{2}+2(\Psi\Theta\tilde{\theta}^{*}\xi,\eta)
+|\eta|^{2}=|\theta\xi|^{2}-|\Psi\Theta\tilde{\theta}^{*}\xi|^{2},
$$
which is even stronger than \eqref{eq3.2.17.1}.

\end{remark}

\begin{remark}                           \label{rem3.2.19.3}
We have seen that 
 Assumptions~\ref{asm3.2.17} and~\ref{asm3.2.16}
 follow from Assumption \ref{assumption 1.15.1}.
 In turn
 Assumptions~\ref{asm3.2.17} and~\ref{asm3.2.16}
in combination with Assumption \ref{asm3.2.15}
imply Assumption \ref{assumption 1.15.1} perhaps
with a different constant  in the latter. 

To show this, we take 
$\zeta=(\xi,\eta)\in
\bR^{d}\times \bR^{d_{1}-d}$ and
observe that 
$$   
2(\tilde{a}\zeta,\zeta)=   
(\theta \theta^{*} \xi,\xi)+2(\Theta \theta^{*}\xi,
\eta)+(\Theta \Theta^{*}\eta,\eta)
$$
$$
=|\theta^{*} \xi|^{2}+2(\Psi\Theta \theta^{*}\xi,
\tilde{\eta})+\varepsilon(\tilde{\eta},\tilde{\eta})+
(1-\varepsilon)(\Theta \Theta^{*}\eta,\eta)
$$
where   $\tilde{\eta}=\Psi^{-1}\eta$,
 and $\varepsilon\in(0,1)$.
By using the inequality $2(\mu,\nu)+\varepsilon|\mu|^{2}\geq-
\varepsilon^{-1}|\nu|^{2}$
we see that
$$
2(\Psi\Theta \theta^{*}\xi,
\tilde{\eta})+\varepsilon(\tilde{\eta},\tilde{\eta})
\geq-\varepsilon^{-1}|\Psi\Theta\tilde{\theta}^{*}\xi|^{2},
$$
 and by taking $N$ such that $\Psi\leq N(\delta^{ij})$,
for which $\Theta\Theta^{*}\geq N^{ -2}(\delta^{ij})$, 
we conclude
$$
2(\tilde{a}\zeta,\zeta)\geq
|\theta^{*} \xi|^{2}-\varepsilon^{-1}|\Psi\Theta\theta^{*}
\xi|^{2}+(1-\varepsilon)N^{ -2}|\eta|^{2}
$$
$$
\geq\delta|\xi|^{2}+(1-\varepsilon^{-1})|\Psi\Theta\theta^{*}
\xi|^{2}+(1-\varepsilon)N^{ -2}|\eta|^{2},
$$
where the last inequality follows from (\ref{eq3.2.17.1}).
Finally, $\Psi\Theta\theta^{*}$ is a bounded function, 
so that,
for a constant $N_{1}$,
$$
2(\tilde{a}\zeta,\zeta)\geq (\delta+N_{1}(1-\varepsilon^{-1}))
|\xi|^{2}+(1-\varepsilon)N^{ -2}|\eta|^{2}.
$$
For $\varepsilon$ sufficiently close to $1$ the last expression
 is greater than $\delta_{1}|\zeta|^{2}$ with a constant
$\delta_{1}>0$, which is equivalent to the
uniform ellipticity of 
$\tilde{a}$.

\end{remark}

 Before stating the next assumption we remind the reader
that, for  $\gamma\in\bR$ and $u\in C^{\infty}_{0}(\bR^{d})$
one introduces $(1-\Delta)^{-\gamma/2}u$ by means
of the Fourier transform. Then, for $p\in(1,\infty)$, one defines
the spaces of Bessel potential $H^{\gamma}_{p}
(\bR^{d})$
as the set of distributions obtained
as the closure of $C^{\infty}_{0}(\bR^{d})$
with respect to the norm 
$$
\|u\|_{H^{\gamma}_{p}(\bR^{d})}:=
\|(1-\Delta)^{ \gamma/2}u\|_{L_{p}(\bR^{d})}.
$$
One important and highly nontrivial piece
of information is that 
$$
H^{1}_{p}(\bR^{d})=W^{1}_{p}(\bR^{d}):=\{u\in L_{p}(\bR^{d}):
\nabla u\in L_{p}(\bR^{d})\}
$$
and
\begin{equation}
                                             \label{1.16.2}
\|u\|_{H^{1}_{p}(\bR^{d})}\sim\|u\|_{W^{1}_{p}(\bR^{d})}
:=\|u\|_{L_{p}(\bR^{d})}+\|\nabla u\|_{L_{p}(\bR^{d})}.
\end{equation}

\begin{assumption}                            \label{asm3.2.18}
The random vectors $x_{0}$ and $y_{0}$ are 
independent of the process
$w_{t}$. The conditional distribution of $x_{0}$ given $y_{0}$ has a
density, which we denote by $\pi_{0}(x)=\pi_{0}(\omega ,x)$. We have
$p\geq2$ and
$\pi_{0} \in L_{p}(\Omega , H_{p}^{1-2/p}(\bR^{d}))$
(actually, we need slightly less, see Remark \ref{remark 1.4.3}).
\end{assumption}

Next we introduce few more notation.
Let
$$
\Psi_{t}=\Psi(t,y_{t}),\quad\Theta_{t}=\Theta(t,y_{t}),\quad
a_{t}(x) =\frac{1}{2}\theta \theta^{*}(t,x,y_{t}),
\quad b_{t}(x)=b(t,x,y_{t}),
$$
$$
\sigma_{t}(x) =\theta(t,x,y_{t}) \Theta^{*}_{t}\Psi_{t},\quad
\beta_{t}(x) =\Psi_{t}B(t,x,y_{t}).
$$
In the remainder of the article we use the notation
$$
D_{i}=\frac{\partial}{\partial x^{i}}
$$
only for $i=1,...,d$ and set
\begin{equation}                             \label{eq3.2.19.2}
L_{t}( x) = a^{ij}_{t}(x)D_{i}D_{j}+
 b^{i}_{t}(x)D_{i}\,,
\end{equation}
$$
L^{*}_{t}( x)u_{t}(x) = 
D_{i}D_{j}(
a^{ij}_{t}( x)u_{t}(x) )
- D_{i}(b^{i}_{t}( x)u_{t}(x) )
$$
\begin{equation}                     \label{e3.2.19.2}
=D_{j}\big(
a^{ij}_{t}(x) D_{i}u_{t}(x)
-b^{j}_{t}(x)u_{t}(x)+u_{t}(x)D_{i}
a^{ij}_{t}(x)\big),
\end{equation}
\begin{equation}                     \label{1.27.8}
\Lambda^{k }_{t}( x)u_{t}(x)  =\beta^{k}_{t}( x)u_{t}(x) +
 \sigma^{ik}_{t}( x) D_{i}u_{t}(x),
\end{equation}
$$
\Lambda^{k*}_{t}( x)u_{t}(x)  =\beta^{k}_{t}( x)u_{t}(x) -
D_{i}(\sigma^{ik}_{t}( x) u_{t}(x) )
$$
\begin{equation}\label{eq3.2.19.3}
=-\sigma^{ik}_{t}(x)D_{i}u_{t}(x)
+(\beta^{k}_{t}(x)-D_{i}
\sigma^{ik}_{t}(x))u_{t}(x),
\end{equation}
where $t\in[0,T]$,   $x\in \bR^{d}$, $k=1,...,d_{1}-d$,
and as above
 we use the summation convention over all ``reasonable''
values of repeated indices,
so that the summation in (\ref{eq3.2.19.2}), (\ref{e3.2.19.2}),
\eqref{1.27.8},
and (\ref{eq3.2.19.3}) is done for $i,j=1,...,d$ (whereas in
(\ref{eq3.2.19.1}) for $i,j=1,...,d_{1}$).
Observe that Lipschitz continuous functions have
bounded generalized derivatives and by
$$
D_{i}
a^{ij}_{t},\quad D_{i}
\sigma^{ik}_{t}
$$
we mean these derivatives. From Remark~\ref{rem3.2.17.2} 
we have that the operator
$L$ defined by \eqref{eq3.2.19.2} is uniformly elliptic 
with constant
of ellipticity $\delta$.

Finally, by $\cF_{t}^{y}$ we denote the completion of 
$\sigma \{ y_{s}:s\leq t\}$
with respect to $P,\cF$.

Let us consider the following initial value problem
\begin{equation}
                                        \label{eq3.2.20}
d\bar{\pi }_{t}(x)=L^{*}_{t}(x )\bar{\pi }_{t}(x)\,dt+
 \Lambda^{k*}_{t}(x )\bar{\pi }_{t}(x)
 \Psi^{kr}_{t}\,dy^{r}_{t},
\end{equation}
$$
\bar{\pi }_{0}(x)= \pi  _{0}(x),
$$
where $t\in[0,T]$, $x\in\bR^{d}$, and
$\bar{\pi }_{t}(x)=\bar{\pi }_{t}(\omega,x)$. 
Equation \eqref{eq3.2.20}
is called the Duncan-Mortensen-Zakai or just the Zakai equation.

We understand this equation and the initial
condition in the following sense.
We are looking for a function $\bar{\pi}=\bar{\pi}_{t}(x)
=\bar{\pi}_{t}(\omega,x)$, $\omega\in\Omega$,
$t\in[0,T]$, $x\in\bR^{d}$, such that

(i) For each $(\omega,t)$, $\bar{\pi}_{t}(\omega,x)$
is a generalized function on $\bR^{d}$,

(ii) We have $\bar{\pi}\in L_{p}(\Omega\times
[0,T],\cP,H^{1}_{p}(\bR^{d}))$,

(iii) For each $\varphi\in C^{\infty}_{0}(\bR^{d})$
with probability one for all $t\in[0,T]$ it holds that
$$
(\bar{\pi}_{t},\varphi)=(\pi_{0},\varphi)
-\int_{0}^{t}(
a^{ij}_{t} D_{i}\bar{\pi}_{t}
-b^{j}_{t}\bar{\pi}_{t}+\bar{\pi}_{t}D_{i}
a^{ij}_{t},D_{j}\varphi)\,dt
$$
\begin{equation}
                                                 \label{1.16.1}
-\int_{0}^{t}
(\sigma^{ik}_{t}D_{i}\bar{\pi}_{t}
+( D_{i}
\sigma^{ik}_{t}-\beta^{k}_{t})\bar{\pi}_{t},\varphi)\Psi^{kr}_{t}
\big(B^{r}(t,z_{t})\,dt+\Theta^{rs}(t,y_{t})\,dw^{s}_{t}\big),
\end{equation}
where by $(f ,\varphi)$ we mean the action
of a generalized function $f$ on $\varphi$, in particular,
if $f$ is a locally summable,
$$
(f,\varphi)=\int_{\bR^{d}}f(x)\varphi(x)\,dx.
$$

Observe that all expressions in \eqref{1.16.1} are well defined
due to the fact that the coefficients 
of $\bar{\pi}$ and of  
$D_{i}\bar{\pi}$
are bounded
and appropriately measurable and $\bar{\pi},
D_{i}\bar{\pi}\in L_{p}(\Omega\times[0,T],
\cP,L_{p}(\bR^{d}))$
(see \eqref{1.16.2}).

Hence, equation \eqref{eq3.2.20} has the same form as \eqref{11.13.1}
and  
the existence and uniqueness part of   Lemma 
\ref{lm3.2.21} below
  follow  from Theorem \ref{theorem 12.7.1} and 
Remark \ref{remark 1.4.3}. The second assertion of the lemma
follows from Theorem \ref{theorem 1.4.3}.

In all what follows in the main 
part of the article we suppose that
   Assumptions \ref{asm3.2.15}, \ref{assumption 1.15.1}, and
  \ref{asm3.2.18}
are satisfied.

\begin{lemma}
                                            \label{lm3.2.21} 
There exists a unique solution $\bar{\pi}$ 
of \eqref{eq3.2.20} with initial
condition $\pi_{0}$
in the sense explained above. In addition,   
$
\bar{\pi}_{t}\geq0
$
for all $t\in[0,T]$ (a.s.).
\end{lemma}

Here is a basic result of filtering theory
for partially observable diffusion
processes. Its relation to the previously
known ones is discussed above.

\begin{theorem}
                                               \label{thm3.2.22}
Let $\bar{\pi }$ be the function from Lemma \ref{lm3.2.21}.
Then 
\begin{equation} 
                                             \label{1.27.07} 
0<\int_{\bR^{d}}\bar{\pi }_{t}(x)\,dx=(\bar{\pi }_{t},1)<\infty 
\end{equation}
 for 
all $t\in[0,T]$  (a.s.)
and for any $t\in[0,T]$ and real-valued,
bounded or nonnegative, (Borel) measurable function $f$ 
given on $\bR^{d}$
\begin{equation} 
                                             \label{10.13.3}
E[f(x_{t})|\cF_{t}^{y}]=
\frac{(\bar{\pi }_{t},f)}
{(\bar{\pi }_{t},1)}\quad
\text{(a.s.).}
\end{equation}
\end{theorem}

Equation \eqref{10.13.3} shows (by definition) that
$$
\pi_{t}(x):=\frac{\bar{\pi }_{t}(x) }
{(\bar{\pi }_{t},1)}
$$
is a conditional density of distribution of
$x_{t}$ given $y_{s},s\leq t$. Since, generally,
$(\bar{\pi }_{t},1)\ne1$,  one calls $\bar{\pi }_{t}$
an unnormalized conditional density of distribution of
$x_{t}$ given $y_{s},s\leq t$.

The following is a direct corollary of 
Theorem \ref{theorem 1.17.1}.

\begin{theorem}
                                        \label{theorem 1.19.1}
Let $\pi_{0}$ be a nonrandom 
function and $\pi_{0}\in H^{1-2/p}_{p}(\bR^{d})$
for all $p\geq2$, which happens
for instance, if $\pi_{0}$ is a
Lipschitz
continuous function with compact support.
Then for any $\varepsilon\in(0,1/2)$ almost surely
$\bar{\pi}_{t}(x)$ is $1/2-\varepsilon$ H\"older 
continuous in $t$ with a constant independent of $x$,
$\bar{\pi}_{t}(x)$ is $1-\varepsilon$ H\"older 
continuous in $x$ with a constant independent of $t$,
and the above mentioned (random) constants have all moments.

\end{theorem}
 
In filtering theory usually the following theorem
is proved before anything else is done. We  do not need 
it 
for proving the above results and  
 give the proof  
just to show
that the   
$L_{p}$-theory of SPDEs allows one to 
get all basic results from filtering theory.

Historically, $P_{t}[\beta]$ was introduced by \eqref{1.28.4}
and shown to have (a modification possessing)
appropriate measurability properties. 
Then $\bar{\pi}_{t}$ used to be defined as the density
of conditional distribution of $x_{t}$ given $\cF^{y}_{t}$
divided by  an appropriate modification of
\begin{equation}
                                      \label{1.29.2}
E(\rho_{t}\mid\cF^{y}_{t}) ,
\end{equation}
  where
$$
\rho_{t}=
\exp(-\int_{0}^{t}
\tilde{\beta}_{s}\,d\tilde{w}_{s}
-\tfrac{1}{2}\int_{0}^{t}|\tilde{\beta}_{s}|^{2}\,ds),
\quad\tilde{\beta}_{s}=\beta_{s}(x_{s}),
\quad  \tilde{w}_{t}=\int_{0}^{t}\Psi_{s}\Theta_{s}\,dw_{s}.
$$
In this case
  $(\bar{\pi}_{t}
,1)^{-1}$ turns out to be this same  appropriate modification of
\eqref{1.29.2}
(cf. our \eqref{1.29.1}).

The most surprising statements in Theorem \ref{theorem 1.28.1}
are assertions (iv) and (v). In (iv) the difference
of two Wiener processes $\check{w}_{t}$ and $\tilde{w}_{t}$
(that the latter is a Wiener process is checked in the proof
of Lemma \ref{lemma 10.21.4})
is asserted to be a differentiable nontrivial function.

Assertion (v) shows that  \eqref{1.29.2}, which is
a conditional expectation of a martingale, is again a martingale
and, moreover, while evaluating it we can just put 
conditional expectations of $\tilde{\beta}_{s}$
given $\cF^{y}_{s}$ in place of $\tilde{\beta}_{s}$
in the expression of $\rho_{t}$ with simultaneous
replacement of $\tilde{w}$ with $\check{w}$.
 
\begin{theorem}
                                     \label{theorem 1.28.1}

(i) The process $( \bar{\pi}_{t},1)$ is  
  continuous in $t$ (a.s.) 
and (a.s.) for all $t\in[0,T]$
\begin{equation}
                                                 \label{1.28.1}
 ( \bar{\pi}_{t},1)=(\pi_{0},1)+
\int_{0}^{t}( \bar{\pi}_{s},\beta^{k}_{s})\Psi^{kr}_{s}
B^{r}(t,z_{s})\,ds+\int_{0}^{t}( \bar{\pi}_{s},\beta^{k}_{s})
\Psi^{kr}_{s}\Theta^{rn}
(t,y_{s})\,dw^{n}_{s}.
\end{equation}

(ii) The process $\bar{\pi}_{t}$ is a continuous
$L_{1}$-valued process (a.s.).

(iii) Introduce $P_{t}[\beta]=(P_{t}[\beta^{1}],...,
P_{t}[\beta^{d_{1}-d}])$ by
$$
P_{t}[\beta]=(\bar{\pi}_{t},1)^{-1}\int_{\bR^{d}}
\beta_{t}(x)\bar{\pi}_{t}(x)\,dx
=(\bar{\pi}_{t},1)^{-1}\Psi(t,y_{t})\int_{\bR^{d}}
B(t,x,y_{t})\bar{\pi}_{t}(x)\,dx.
$$
Then $P_{t}[\beta]$ is a jointly measurable bounded
$\cF^{y}_{t}$-adapted process on $[0,t]$ (a.s.) and
for each $t\in[0,T]$
\begin{equation}
                                         \label{1.28.4}
P_{t}[\beta]=E(\beta_{t}(x_{t})\mid \cF^{y}_{t})
\quad\text{(a.s.)}.
\end{equation}

(iv) The process
$$
\check{w}_{t}=\tilde{w}_{t}+\int_{0}^{t}
(\beta_{s}(x_{s})-P_{s}[\beta])\,ds
$$
is a $(d_{1}-d)$-dimensional
Wiener process with respect to  $\cF^{y}_{t}$  
(the so-called innovation process),
where
$$
\tilde{w}_{t}=\int_{0}^{t}\Psi_{s}\Theta_{s}\,dw_{s}.
$$

(v) We have (a.s.) for all $t\in[0,T]$
\begin{equation}
                                         \label{1.28.3}
(\bar{\pi},1)=\exp\big(\int_{0}^{t}
P_{s}[\beta]\,d\check{w}_{s}+\tfrac{1}{2}
\int_{0}^{t}|P_{s}[\beta]|^{2}\,ds\big),
\end{equation}
so that
$$
(\bar{\pi},1)^{-1}=\exp\big(-\int_{0}^{t}
P_{s}[\beta]\,d\check{w}_{s}-\tfrac{1}{2}
\int_{0}^{t}|P_{s}[\beta]|^{2}\,ds\big)
$$
is an exponential martingale,
and for each $m>0$
\begin{equation}
                                         \label{1.28.03}
E\sup_{t\leq T}(\bar{\pi},1)^{m}+
E\sup_{t\leq T}(\bar{\pi},1)^{-m}<\infty.
\end{equation}

\end{theorem}

 \mysection{Proof of Theorem \protect\ref{thm3.2.22}}
                                           \label{section 1.31.1}

We will use some notion and results from the 
theory of SPDEs, which are recalled  
in Section \ref{section 1.31.3}.
From now on we drop $\bR^{d}$ in notation
like $H^{\gamma}_{p}(\bR^{d})$ and $L_{p}(\bR^{d})$.

\begin{remark}
                                     \label{remark 1.4.3}

The assumption that $\pi_{0} \in L_{p}(\Omega , H_{p}^{1-2/p})$
is only needed to guarantee (see the proof of Theorem 5.1
of \cite{K99}) that there exists a $\psi\in\cH^{1}_{p}(T)$
such that $\psi_{0}=\pi_{0}$,
$$
d\psi_{t}=\Delta\psi_{t}\,dt=D_{i}
f^{i}_{t}\,dt,\quad (f^{i}_{t}=D_{i}
\psi_{t}),
$$
$$
\|\psi\|_{\bH^{1}_{p}(T)}^{p}\leq NE\|\pi_{0}\|_{H^{1-2/p}
_{p}}^{p}
$$
with $N$ independent of $\pi_{0}$.
\end{remark}

As   
is mentioned before Lemma \ref{lm3.2.21},
by Theorem \ref{theorem 12.7.1} and 
Remark \ref{remark 1.4.3},
there exists a unique solution  
$\bar{\pi}\in\cH^{1}_{p}(T)$ 
of \eqref{eq3.2.20} with initial
condition $\pi_{0}$. 
By Theorem \ref{theorem 1.4.3},
$
\bar{\pi}_{t}\geq0
$
for all $t\in[0,T]$ (a.s.). By Theorem
\ref{theorem 1.17.1},  $\bar{\pi}_{t}$
is a continuous $L_{p}$-valued process and
\begin{equation}
                                                 \label{1.27.11}
E\sup_{t\in[0,T]}\|\bar{\pi}_{t}\|^{p}_{L_{p}}\,dt<\infty.
\end{equation}

Now, we prove three auxiliary results.
\begin{lemma}
                                        \label{lemma 1.6.1}
Let $\xi_{t},\xi^{n}_{t}$, $n=1,2,...$, $t\in[0,T]$,
be $k$-dimensional continuous semimartingales
 such that, for any $t\in[0,T]$,
$\xi^{n}_{t}\to\xi_{t}$ in probability as $n\to\infty$.
Assume that
$$
\xi^{n}_{t}=\xi^{n}_{0}+\int_{0}^{t}\alpha^{n}_{s}\,ds
+m^{n}_{t},\quad
\xi_{t}=\xi_{0}+\int_{0}^{t}\alpha _{s}\,ds
+m_{t},
$$
where $\alpha_{t}$ and $\alpha^{n}_{t}$ are predictable processes
bounded by the same nonrandom constant and $m_{t}$ and $m^{n}_{t}$
are martingales such that
$$
\langle m^{ni},m^{nj}\rangle_{t}=\int_{0}^{t}
\gamma^{nij}_{s}\,ds,\quad
\langle m^{i},m^{j}\rangle_{t}=\int_{0}^{t}
\gamma^{ij}_{s}\,ds,\quad i,j=1,...,k,
$$
where $\gamma^{n}_{t}:=(\gamma^{nij}_{t})$
 and $\gamma_{t}:=(\gamma^{ij}_{t})$ are predictable 
matrix-valued processes
bounded by the same nonrandom constant and such that
$(\gamma^{n}_{t})^{-1}$ and $(\gamma_{t})^{-1}$
exist and are also bounded by the same nonrandom constant.

Assume that on $[0,T]\times\bR^{l}\times\bR^{k}$ we are given
functions $f^{n}_{t}(x,y)$ and $f_{t}(x,y)$ such that
they are uniformly bounded and $f^{n}\to f$ in measure
as $n\to\infty$.

Then $f^{n}_{t}(x,\xi^{n}_{t})\to
f_{t}(x,\xi_{t})$ in measure on $\Omega\times[0,T]
\times\bR^{l}$.

\end{lemma}

Proof. It suffices to show that any subsequence $\{n'\}$
of integers has a subsequence $\{n''\}$ such that
$f^{n''}_{t}(x,\xi^{n''}_{t})\to
f_{t}(x,\xi_{t})$ in measure.  Since any subsequence
$\{n'\}$ has a subsequence $\{n''\}$ such that
$f^{n''}\to f$ almost everywhere,  by 
having in mind renumbering if needed,
we may assume that for the original sequence we have
$f^{n}\to f$ almost everywhere. In that case
for almost any $x\in\bR^{l}$, $f^{n}_{t}(x,y)\to f_{t}(x,y)$
and, if we prove that for each such $x$ we have
$f^{n}_{t}(x,\xi^{n}_{t})\to
f_{t}(x,\xi_{t})$ in measure on $\Omega\times[0,T]$, then
$$
E\int_{0}^{T}|f^{n}_{t}(x,\xi^{n}_{t})
-f_{t}(x,\xi_{t})|\,dt\to0,
$$
which after being integrated with respect to $x$ would shows
that $f^{n}_{t}(x,\xi^{n}_{t})\to
f_{t}(x,\xi_{t})$ in measure on $\Omega\times[0,T]
\times\bR^{l}$.

It follows that we only need to prove that,
if  on $[0,T]\times\bR^{k}$ we are given
functions $f^{n}_{t}(y)$ and $f_{t}(y)$ such that
they are uniformly bounded and $f^{n}\to f$ 
$(t,y)$-almost everywhere
as $n\to\infty$, then 
\begin{equation}
                                                 \label{1.6.1}
E\int_{0}^{T}|f^{n}_{t}(\xi^{n}_{t})
-f_{t}(\xi_{t})|\,dt\to0.
\end{equation}
Furthermore, since the coefficients $\alpha^{n}$, $\alpha$,
$\gamma^{n}$, and $\gamma$ are uniformly bounded
$$
\sup_{n}\sup_{t\in[0,T]}P(|\xi^{n}_{t}|+|\xi_{t}|\geq R)
\leq R^{-2}\sup_{n}\sup_{t\in[0,T]}E(|\xi^{n}_{t}|^{2}
+|\xi_{t}|^{2})
\to0
$$
as $R\to\infty$. Therefore, if for any $R\in(0,\infty)$
we know that \eqref{1.6.1} is true
provided that $f^{n}_{t}(y)$ and $f_{t}(y)$ vanish
for $|y|\geq R$, then by applying this result in the general case
to  $f^{n}_{t}(y)I_{|y|<R}$ and $f_{t}(y)I_{|y|<R}$
we would obtain that
$$
\nlimsup_{n\to\infty} E\int_{0}^{T}|f^{n}_{t}(\xi^{n}_{t})
-f_{t}(\xi_{t})|\,dt\leq NR^{-2},
$$
where $N$ is independent of $R$. This would imply
\eqref{1.6.1} in the general case. This shows that
without restricting generality we may assume that
for an $R\in(0,\infty)$ the functions
$f^{n}_{t}(y)$ and $f_{t}(y)$ vanish
if $|y|\geq R$.

Now observe that the left-hand side of \eqref{1.6.1}
is majorated by $I_{n}+J_{n}$, where
$$
I_{n}=E\int_{0}^{T}|f^{n}_{t}(\xi^{n}_{t})
-f_{t}(\xi^{n}_{t})|\,dt,\quad J_{n}
=E\int_{0}^{T}|f_{t}(\xi^{n}_{t})-f_{t}(\xi_{t})|\,dt.
$$

We recall a result of \cite{K77} implying that
for any $g\in L_{k+1}([0,T]\times\bR^{k})$ we have
$$
E\int_{0}^{T}(|g_{t}(\xi^{n}_{t})|+
|g_{t}(\xi_{t})|)\,dt\leq N\|g\|_{L_{k+1}([0,T]\times\bR^{k})},
$$
where $N$ is independent of $n$ and $g$.
We apply this result to $g=f^{n}-f$ 
and observe that these functions are uniformly bounded,
 vanish for $|y|\geq R$, and tend to zero in measure.
Hence, their $L_{k+1}([0,T]\times\bR^{k})$-norms tend to zero.
This implies that $I_{n}\to0$.

Next, notice that for any function $g$
$$
J_{n}\leq
E\int_{0}^{T}|g_{t}(\xi^{n}_{t})-g_{t}(\xi_{t})|\,dt
$$
$$
+E\int_{0}^{T}|f_{t}(\xi^{n}_{t})-g_{t}(\xi^{n}_{t})|\,dt
+E\int_{0}^{T}|f_{t}(\xi_{t})-g_{t}(\xi_{t})|\,dt
$$
implying that
\begin{equation}
                                                 \label{1.6.2}
\nlimsup_{n\to\infty}J_{n}\leq
\nlimsup_{n\to\infty}E\int_{0}^{T}|g_{t}
(\xi^{n}_{t})-g_{t}(\xi_{t})|\,dt
+N\|f-g\|_{L_{k+1}([0,T]\times\bR^{k})},
\end{equation}
where $N$ is independent of $g$. For any $\varepsilon>0$
we can find a smooth $g$ such that the second term on the right
in \eqref{1.6.2}
will be less than $\varepsilon$. In addition, the first term
vanishes for smooth $g$ since $\xi^{n}_{t}
\to\xi_{t}$ in probability for any $t$. Since $\varepsilon$
is arbitrary, it follows that the left-hand side of
\eqref{1.6.2} equals zero. The lemma is proved.

The following result with its proof is an adaptation
of Lemma 5.1 of \cite{KZ} and its proof.
\begin{lemma}
                                          \label{lemma 10.21.4}
The function $\bar{\pi}_{t}$ 
  is $\cF^{y}_{t}$-adapted.
\end{lemma}

Proof. Define 
$$
\tilde{\beta}_{t}=\beta_{t}(x_{t})=
\Psi_{t}B(t,z_{t}),\quad
\hat{w}_{t}=\int_{0}^{t}\Psi_{s}\,dy_{s},\quad
\tilde{w}_{t}=\int_{0}^{t}\Psi_{s}\Theta_{s}\,dw_{s}.
$$
Since $\Psi_{t}$ is $\cF^{y}_{t}$-adapted, the process
$\hat{w}_{t}$ is $\cF^{y}_{t}$-adapted too. Furthermore,
$\Psi_{s}\Theta_{s}\Theta^{*}_{s}\Psi_{s}$ is a unit matrix
so that by L\'evy's theorem $\tilde{w}_{t}$ is a Wiener
process.
We want to change the probability measure so that
$\hat{w}_{t}$ would become a Wiener process with respect
to this new measure. Define
\begin{equation}
                                           \label{1.27.1}
\rho_{t}=
\exp(-\int_{0}^{t}
\tilde{\beta}_{s}\,d\tilde{w}_{s}
-\tfrac{1}{2}\int_{0}^{t}|\tilde{\beta}_{s}|^{2}\,ds),\quad
Q(d\omega)=\rho_{T}(\omega)\,P(d\omega).
\end{equation}
The process $\rho_{t}$ is an exponential local martingale.
Since $\tilde{\beta}$ is bounded,  $\rho_{t}$ is
square integrable, so that $Q$ is a probability
measure. Since
$$
d\hat{w}_{t}=\tilde{\beta}_{t}\,dt+d\tilde{w}_{t} 
$$
and $\tilde{w}_{t}$ is a Wiener process on $(\Omega,\cF,P)$,
by Girsanov's theorem, $\hat{w}_{t}$, $t\in[0,T]$,
 is a Wiener process
on $(\Omega,\cF,Q)$ with respect to the filtration
$\{\cF_{t}\}$. As has been noticed before, it is 
$\cF^{y}_{t}$-adapted and, obviously, 
$$
\cF^{y}_{t}
\subset\cF_{t},
$$
 so that $(\hat{w}_{t},\cF^{y}_{t})$
is a Wiener process. Now rewrite (\ref{eq3.2.20}) as
\begin{equation}
                                           \label{10.22.1}
d\bar{\pi }_{t}(x)=L^{*}_{t}(x )\bar{\pi }_{t}(x)\,dt+
 \Lambda^{k*}_{t}(x )\bar{\pi }_{t}(x)  \,d\hat{w}^{k}_{t},
\end{equation}
and consider this equation relative to
$(\Omega,\cF,\cF^{y}_{t},Q)$. 

By Theorem \ref{theorem 12.7.1} and Remark \ref{remark 1.4.3} 
equation \ref{10.22.1} with initial data $\pi_{0}$
has a unique   $\cF^{y}_{t}$-adapted solution
belonging to $\cH^{1}_{p}(
\cF^{y}_{\cdot},Q,T)\subset\cH^{1}_{p}(
\cF_{\cdot},Q,T)$, where by $\cH^{1}_{p}(
\cF^{y}_{\cdot},Q,T)$
we mean the space $\cH^{1}_{p}(T)$ constructed on the basis
of the new probability measure $Q$ and filtration
$\cF^{y}_{\cdot}$. We denote by
$\tilde{\pi}_{t}$ this  solution.

We have already mentioned that $\bar{\pi}
\in\cH^{1}_{p}(\cF_{\cdot},P,T)$.
We want to derive that $\bar{\pi}_{t}$ is $\cF^{y}_{t}$-adapted
from the uniqueness by showing that $\bar{\pi}=\tilde{\pi}$
because both are $\cF_{t}$-adapted solutions of the same equation.
The only obstacle is that the norms in 
$\cH^{1}_{p}(\cF_{\cdot},Q,T)$ and $\cH^{1}_{p}( T)$ are different.
To overcome this obstacle,
  we are going to use stopping times.

For integers $n$ define 
$$
\tau (n)=T\wedge\inf\{t\geq0:\int_{0}^{t} 
\|\tilde{\pi}_{s }\|^{p}_{H^{1}_{p}} 
\,ds\geq n\}.
$$
Obviously, $\tau(n)$ are  $\cF^{y}_{t}$-stopping
times and $\cF_{t}$-stopping
times. Furthermore,
$$
\|\tilde{\pi}\|_{\bH^{1}_{p}(\cF_{\cdot},P,\tau(n))}^{p}=
E \int_{0}^{\tau(n)}\|\tilde{\pi}_{s }\|^{p}_{H^{1}_{p}} 
\,ds\leq n<\infty.
$$
This and the equation (cf.~\eqref{10.22.1})
$$
d\tilde{\pi }_{t}(x)=\big[L^{*}_{t}(x )\tilde{\pi }_{t}(x)
+\tilde{\beta}^{k}_{t}
\Lambda^{k}_{t}(x )\tilde{\pi }_{t}(x)\big]\,dt+
 \Lambda^{k}_{t}(x )\tilde{\pi }_{t}(x)  \,d\tilde{w}^{k}_{t}
$$
show that,
 $\tilde{\pi}\in\cH^{1}_{p}(\cF_{\cdot},P,\tau(n))$.
By the above mentioned uniqueness, $\tilde{\pi}_{t}=\bar{\pi}_{t}$
on $\opar0,\tau(n)\cbrk$ (a.e.). Since both functions are 
continuous in $t\in[0,T]$ (Theorem \ref{theorem 1.17.1} (i)),
we have that 
$$
\tilde{\pi}_{t}I_{0<t\leq\tau(n)}\quad\text{and}\quad
 \bar{\pi}_{t}I_{0<t\leq\tau(n)}
$$
 are indistinguishable,
and since one of them is $\cF^{y}_{t}$-adapted, so is the other.
We conclude that $\bar{\pi}_{t}I_{0<t\leq\tau(n)}$
is $\cF^{y}_{t}$-adapted, which after letting $n\to\infty$
yields the result. The lemma is proved.

Assertion  of the following lemma is a very particular
case of one of the assertions of Theorem \ref{theorem 1.28.1}.
 Before stating the lemma
we recall that
$\bar{\pi}_{t}\geq0$ for all $t\in[0,T]$ (a.s.),
so that $( \bar{\pi}_{t},1)$ is well defined
(and may be infinite).
\begin{lemma}
                                       \label{lemma 1.27.2}
We have
\begin{equation}
                                             \label{1.27.9} 
E\sup_{t\in[0,T]}(\bar{\pi}_{t},1)^{1/2}<\infty.
\end{equation}
 \end{lemma}

Proof.  
For $\varphi\in C^{\infty}_{0}(\bR^{d})$
one can rewrite \eqref{1.16.1} as
$$
(\bar{\pi}_{t},\varphi)=(\pi_{0},\varphi)
+\int_{0}^{t}(
 \bar{\pi}_{s},L_{s}\varphi)\,ds
$$
\begin{equation}
                                                 \label{1.27.10}
+\int_{0}^{t}
( \bar{\pi}_{s},\Lambda^{k}_{s}\varphi)\Psi^{kr}_{s}
\big(B^{r}(s,z_{s})\,ds+\Theta^{rn}(s,y_{s})\,dw^{n}_{s}\big).
\end{equation}
Using \eqref{1.27.11} and  an obvious
passage to the limit, it is easy to prove that 
\eqref{1.27.10} holds not only for 
$\varphi\in C^{\infty}_{0}(\bR^{d})$, but also
for $\varphi\in W^{2}_{q}$ with $q=p/(p-1)$.

On $\bR^{d}$ for $m=1,2,...$ introduce the functions
$$
\varphi(x)=(1+|x|^{2})^{-d},\quad\varphi_{m}(x)=\varphi(x/m).
$$
Observe that for a constant $N$ it holds that
\begin{equation}
                                                 \label{1.28.6}
|D_{i}\varphi_{m}|+ m |D_{i}D_{j}\varphi_{m}|
\leq Nm^{-1}\varphi_{m}
\end{equation}
on $\bR^{d}$ for all $m$. In particular,
\begin{equation}
                                                 \label{1.27.12}
2L_{t}\varphi_{m}\leq N_{0}\varphi_{m},\quad
2|\Psi^{kr}_{t}B^{r}(t,z_{t})\Lambda^{k}_{t}\varphi_{m}|\leq
N_{0}\varphi_{m},
\end{equation}
where $N_{0}$ is a constant independent of $m$ and the arguments
of the functions involved.

By plugging in \eqref{1.27.10} the function $\varphi_{m}$
in place of $\varphi$, 
we obtain
$$
(\bar{\pi}_{t},\varphi_{m})=(\pi_{0},\varphi_{m})
+\int_{0}^{t}(
 \bar{\pi}_{s},L_{s}\varphi_{m})\,ds
$$
\begin{equation}
                                                 \label{1.28.5}
+\int_{0}^{t}
( \bar{\pi}_{s},\Lambda^{k}_{s}\varphi_{m})\Psi^{kr}_{s}
\big(B^{r}(s,z_{s})\,ds+\Theta^{rn}(s,y_{s})\,dw^{n}_{s}\big).
\end{equation}
By 
  using It\^o's formula
for transforming 
\begin{equation}
                                                 \label{1.27.13}
(\bar{\pi}_{t},\varphi_{m})e^{-N_{0}t},
\end{equation}
and   using \eqref{1.27.12} we 
see that
$$
d\big[(\bar{\pi}_{t},\varphi_{m})e^{-N_{0}t}\big]=
e^{-N_{0}t}( \bar{\pi}_{t},\Lambda^{k}_{t}\varphi_{m})\Psi^{kr}_{t}
 \Theta^{rn}(t,y_{t})\,dw^{n}_{t} 
$$
$$
+e^{-N_{0}t}[(
 \bar{\pi}_{t},L_{t}\varphi_{m})+
( \bar{\pi}_{s},\Lambda^{k}_{s}\varphi_{m})\Psi^{kr}_{s}
 B^{r}(s,z_{s})-N_{0}(\bar{\pi}_{s},\varphi_{m})]\,dt
$$
$$
\leq e^{-N_{0}t}( \bar{\pi}_{t},\Lambda^{k}_{t}
\varphi_{m})\Psi^{kr}_{t}
 \Theta^{rn}(t,y_{t})\,dw^{n}_{t} .
$$
It follows that   process
\eqref{1.27.13} is a supermartingale. It is continuous
and nonnegative. Therefore,
$$
E\sup_{t\in[0,T]}e^{-N_{0}t}
\big(\int_{\bR^{d}}\varphi_{m}\bar{\pi}_{t}(x)
\,dx\big)^{1/2}\leq 2\big(E\int_{\bR^{d}}\varphi_{m}
\bar{\pi}_{0}(x)\,dx\big)^{1/2}\leq 2.
$$
 Upon letting
$m\to\infty$ and using the monotone
convergence theorem we come to \eqref{1.27.9} and the lemma 
is proved.

{\bf Proof of Theorem \ref{thm3.2.22}}. Take
a nonnegative $\zeta\in C^{\infty}_{0}(\bR^{d_{1}})$, which
integrates to one and for $n=1,2,...$ set
$$
\zeta_{n}(z)=n^{d_{1}}\zeta(nz).
$$
Also introduce mollifications of one of the coefficients of 
\eqref{eq3.2.14} by
$$
\theta^{(n)}(t,z) =\zeta_{n}(z)*\theta (t,z),
$$
where the convolutions is taken with respect to $z$.

The function $\zeta$ can be considered as the density
of a random variable.
If needed, we extend our initial probability space
in such a way that it would allow us to introduce
a new  random $\bR^{d_{1}}$-valued vector $\xi$ 
 having density $\zeta$ and such that $\xi$ is independent
of $z_{0}$ and the process  
$w_{t}$, $t\geq0$.

After that, for $n=1,2,...$, we consider the following modification
of \eqref{eq3.2.14}:

\begin{equation}\begin{split}         
                                            \label{1.25.1} 
& dx^{(n)}_{t}=b (t,z^{(n)}_{t}) dt+
\theta ^{(n)}(t,z^{(n)}_{t})dw_{t} \\ 
& dy^{(n)}_{t}=B (t,z^{(n)}_{t}) dt+\Theta (t,y^{(n)}_{t})dw_{t}
\end{split}\end{equation}
with  initial data $x^{(n)}_{0}=x_{0}+n^{-1}\xi$,
$y^{(n)}_{0}=y_{0}$ and $z^{(n)}_{t}=(x^{(n)}_{t},y^{(n)}_{t})$.
Observe that the conditional distribution of $x^{(n)}_{0}$
given $y_{0}$ has a density equal to
$$
\pi^{(n)}_{0} =\zeta_{n}*\pi_{0}.
$$

Since $\theta(t,x,y)$ is Lipschitz in $x$ (even in $(x,y)$)
we have $|\theta(t,z)-\theta^{(n)}(t,z)|\leq Nn^{-1}$,
where $N$ is independent of $n,t,z$. This shows that
system \eqref{1.25.1} satisfies Assumption 
\ref{assumption 1.15.1} for all large $n$.
In addition $\theta^{(n)}$ possesses
enough smoothness in order for the results of
\cite{KZ} to be applicable. For all large $n$, it follows
that, for any  smooth bounded and nonnegative function 
 $c_{t}(y)$ on $[0,T]\times\bR^{d_{1}-d}$
and any  $\varphi\in C^{\infty}_{0}(\bR^{d_{1}})$,
$$
E \varphi(z^{(n)}_{T})\exp(-\int_{0}^{T}c_{s}(y^{(n)}_{s})\,ds)
$$
\begin{equation}
                                              \label{10.21.5}
=E  \rho^{(n)}_{T}
\int_{\bR^{d}}
\varphi(x,y^{(n)}_{T})\bar{\pi}^{(n)}_{T}(x)\,dx
\,\exp(-\int_{0}^{T}c_{s}(y^{(n)}_{s})\,ds),
\end{equation}
where
$\bar{\pi}^{(n)}_{t}$ is the solution of equation 
\eqref{eq3.2.20} corresponding to system \eqref{1.25.1} 
with initial condition $\bar{\pi}^{(n)}_{0}=\pi^{(n)}_{0}$
and $\rho^{(n)}_{t}$ is introduced as in \eqref{1.27.1}
on the basis of \eqref{1.25.1}:
$$
\rho^{(n)}_{t}=\exp(-\int_{0}^{t}
\tilde{\beta}^{(n)}_{s}\,d\tilde{w}^{(n)}_{s}
-\tfrac{1}{2}\int_{0}^{t}|\tilde{\beta}^{(n)}_{s}|^{2}\,ds),
$$
$$
\tilde{w}^{(n)}_{t}=\int_{0}^{t}
\Psi^{(n)}_{s}\Theta^{(n)}_{s}\,dw_{s},\quad
\tilde{\beta}^{(n)}_{t}=\beta ^{(n)}_{t}(x ^{(n)}_{t})
,\quad
 \beta ^{(n)}_{t}(x)=\Psi^{(n)}_{t}B(t,x,y^{(n)}_{t}),
$$
$$
\Theta^{(n)}_{t}=\Theta(t,y^{(n)}_{t}),\quad\Psi^{(n)}_{t}
=\Psi(t,y^{(n)}_{t}).
$$

Later on we will also use the following notation
for other coefficients of equation 
\eqref{eq3.2.20} corresponding to system \eqref{1.25.1}.
Introduce 
$$
a^{(n)}_{t}(x)=\frac{1}{2}\theta^{(n)}
\theta^{(n)*}(t,x,y^{(n)}_{t}),
\quad b^{(n)}_{t}(x)=b(t,x,y^{(n)}_{t}),
$$
$$
\sigma^{(n)}_{t}(x)=\theta^{(n)}(t,x,y^{(n)}_{t})
\Theta^{(n)*}_{t}\Psi^{(n)}_{t}.
$$

Since we know that $\bar{\pi}^{(n)}_{t}\geq0$, it follows from
the validity of
\eqref{10.21.5} for all $\varphi\in C^{\infty}_{0}(\bR^{d_{1}})$,
that it is also valid for all Borel nonnegative 
or bounded $\varphi$.
In particular, for any   $f\in C^{\infty}_{0}(\bR^{d})$  
(independent of $y$) we have
$$
E 
f(x^{(n)}_{T})
 \exp(-\int_{0}^{T}c_{s}(y^{(n)}_{s})\,ds)
$$
\begin{equation}
                                              \label{1.27.2}
=E  \rho^{(n)}_{T}
\int_{\bR^{d}}
f(x )\bar{\pi}^{(n)}_{T}(x)\,dx
\,\exp(-\int_{0}^{T}c_{s}(y^{(n)}_{s})\,ds).
\end{equation}

Our next step is to pass to the limit in \eqref{1.27.2}
as $n\to\infty$. It is a standard fact that for any $m>0$
\begin{equation}
                                              \label{1.27.3}
\lim_{n\to\infty}E\sup_{t\leq T}|z^{(n)}_{t}-z_{t}|^{m}=0,
\end{equation}
which, in particular, implies that the left-hand sides
of \eqref{1.27.2} tend to
$$
E  
f(x _{T})\exp(-\int_{0}^{T}c_{s}(y _{s})\,ds).
$$

Furthermore, the process $\rho^{(n)}_{t}$ is the solution
of the linear equation
$$
d\rho^{(n)}_{t}=-\rho^{(n)}_{t}\gamma^{(n)}_{t}\,dw_{t},
$$
with initial condition $\rho^{(n)}_{0}=1$,
where
$$
\gamma^{(n)}_{t}=\Psi(t,y^{(n)}_{t})B(t,z^{(n)}_{t})
\Psi(t,y^{(n)}_{t})\Theta(t,y^{(n)}_{t}) .
$$
Also introduce 
$$
\gamma _{t}=\Psi(t,y _{t})B(t,z _{t})
\Psi(t,y _{t})\Theta(t,y _{t}) 
$$
and observe that the processes $\gamma^{(n)}_{t}$ and 
$\gamma _{t}$ are bounded.

Furthermore, it follows from \eqref{1.27.3} that for any $m>0$
$$
\lim_{n\to\infty}E\sup_{t\leq T}
|\gamma^{(n)}_{t}-\gamma_{t}|^{m}=0,
$$
which in turn implies that
$$
\lim_{n\to\infty}E\sup_{t\leq T}
|\rho^{(n)}_{t}-\rho_{t}|^{m}=0,
$$
where $\rho_{t}$ is the solution of the equation
$d\rho _{t}=-\rho _{t}\gamma _{t}\,dw_{t}$ with
initial condition $\rho_{0}=1$ and is given
in \eqref{1.27.1}.

To investigate the limit of the remaining factor
on the right in \eqref{1.27.2} we will use
Theorem \ref{theorem 1.4.1}.
By the well-known properties of convolutions
$$
\|\pi^{(n)}_{0} \|^{p}_{H^{1-2/p}}
\leq
\| \pi_{0}\|^{p}_{H^{1-2/p}},\quad
\lim_{n\to\infty}
E\|\pi^{(n)}_{0}-\pi_{0}\|^{p}_{H^{1-2/p}}=0.
$$
This and   Remark \ref{remark 1.4.3} show that
the assumption of Theorem \ref{theorem 1.4.1}
regarding the convergence of the 
initial data for $\bar{\pi}^{(n)}_{t}$ and 
$ \bar{\pi}_{t}$ is satisfied. Furthermore,
there are no free terms
in filtering equations. Therefore, it only remains to check
the appropriate convergence of the coefficients. 
Theorem \ref{theorem 1.4.1} requires the following
convergences in measure $P(d\omega)dtdx$ to hold
on $\Omega\times[0,T]\times\bR^{d}$:
$$
a^{(n)}_{t}(x)\to a_{t}(x),\quad b^{(n)}_{t}(x)\to b_{t}(x)
,\quad D_{i}a^{(n)ij}_{t}(x)\to D_{i}a^{ij}_{t}(x),
$$
$$
\sigma^{(n)}_{t}(x)\to\sigma _{t}(x),
\quad \beta^{(n)}_{t}(x)\to \beta _{t}(x),
\quad D_{i}\sigma^{(n)ik}_{t}(x)\to
D_{i}\sigma^{ ik}_{t}(x).
$$
Relation \eqref{1.27.3} and the assumption that the coefficients
of system \eqref{eq3.2.14}  are Lipschitz
continuous show that, actually, apart from cases involving the
derivatives of $a$ and $\sigma$
all the remaining convergences
hold uniformly in $(t,x)$ almost surely. It is easy to see that
in order to  take care of the terms with derivatives it suffices
to check that
\begin{equation}
                                              \label{1.27.4}
D_{i}\theta^{(n)}(t,x,y^{(n)}_{t})\to
D_{i}\theta (t,x,y _{t})
\end{equation}
in measure for any $i=1,...,d$. Observe that
by the well known properties of
convolutions
$$
D_{i}\theta^{(n)}(t,x,y )\to D_{i}\theta (t,x,y)
$$
for almost all $(t,x,y)$. Therefore, applying Lemma 
\ref{lemma 1.6.1} shows that \eqref{1.27.4} holds.

Now by Theorem \ref{theorem 1.4.1} 
and H\"older's inequality
we conclude
\begin{equation}
                                              \label{1.27.5}
 \lim_{n\to\infty}E\big|\int_{\bR^{d}}
f(x )\bar{\pi}^{(n)}_{T}(x)\,dx-
\int_{\bR^{d}}
f(x )\bar{\pi} _{T}(x)\,dx\big|^{p}=0.
\end{equation}
This along with the above investigation of other
terms in \eqref{1.27.2} yields after letting $n\to\infty$ that
$$
E  
f(x _{T})\exp(-\int_{0}^{T}c_{s}(y _{s})\,ds)=
E  \rho _{T}
(
 \bar{\pi} _{T},f) 
\exp(-\int_{0}^{T}c_{s}(y _{s})\,ds).
$$
The arbitrariness of $c$ leads to
$$
E\big(f(x _{T})\mid \cF^{y}_{T}\big)
=E\big(\rho _{T}
(
 \bar{\pi} _{T},f) \mid \cF^{y}_{T}\big),\quad\text{(a.s.)},
$$
which combined with the $\cF^{y}_{T}$-measurability
of $\bar{\pi} _{T}$ (Lemma \ref{lemma 10.21.4}) shows that
\begin{equation}
                                               \label{1.27.6}
E\big(f(x _{T})\mid \cF^{y}_{T}\big)
=
(
 \bar{\pi} _{T},f) 
E\big(\rho _{T}\mid \cF^{y}_{T}\big) \quad\text{(a.s.)}.
\end{equation}

Observe that on the set of $\omega$ where
\begin{equation}
                                               \label{1.27.7}
E\big(\rho _{T}\mid \cF^{y}_{T}\big)=0
\end{equation}
we have (a.s.)
$$
E\big(f(x _{T})\mid \cF^{y}_{T}\big)=0.
$$
The arbitrariness of $f$ shows that on the said set (a.s.)
$$
1=E\big(1\mid \cF^{y}_{T}\big)=0
$$
and consequently \eqref{1.27.7} can only happen with
probability zero.

Furthermore, by Theorem \ref{theorem 1.4.3} we have
$\bar{\pi}_{t}\geq0$. A standard measure-theoretic
argument then shows that \eqref{1.27.6} holds for all
nonnegative Borel $f$ rather than only for $f\in C^{\infty}_{0}
(\bR^{d})$. By taking $f\equiv1$ we see that
$$
1=
(
 \bar{\pi} _{T},1) 
E\big(\rho _{T}\mid \cF^{y}_{T}\big) \quad\text{(a.s.)}
$$
implying that
\begin{equation}
                                  \label{1.29.1}
\infty>(
 \bar{\pi} _{T},1) >0,\quad
E\big(\rho _{T}\mid \cF^{y}_{T}\big)
= (
 \bar{\pi} _{T},1) ^{-1} \quad\text{(a.s.)}.
\end{equation}
Coming back to  \eqref{1.27.6} we conclude
$$
E[f(x_{T})|\cF_{T}^{y}]=
\frac{ (
 \bar{\pi} _{T},f) }
{(
 \bar{\pi} _{T},1) }\quad
\text{ (a.s.) }
$$
for any nonnegative and any bounded Borel $f$ as well.
Obviously, one can replace here $T$ with any $t\in[0,T]$
and to prove Theorem \ref{thm3.2.22} it only remains
to show that  (a.s.) relation 
\eqref{1.27.07} holds for all $t\in[0,T]$.

The second inequality in \eqref{1.27.07}
holds due to Lemma \ref{lemma 1.27.2}. To prove the first
one it only remains to observe that by the above
for each particular $t\in[0,T]$ with probability one
$$
\int_{\bR^{d}}\bar{\pi}^{p}_{t}(x)\,dx>0
$$
and by Theorem \ref{theorem 1.17.1} the above integral
is  continuous in $t$ with probability one.
The theorem is proved.

\mysection{Proof of Theorem \protect\ref{theorem 1.28.1}}
                                           \label{section 1.31.2}
 
To prove (i) we first show that the right-hand sides of
\eqref{1.28.5}
converge as $n\to\infty$
uniformly in $t\in[0,T]$ in probability
to the right-hand side of \eqref{1.28.1}.
Owing to \eqref{1.28.6} and \eqref{1.27.9}
$$
\int_{0}^{T}|(
 \bar{\pi}_{s},L_{s}\varphi_{m})|\,ds
\leq NTm^{-1}\sup_{s\in[0,T]}(
 \bar{\pi}_{s},1)\to 0\quad\text{(a.s.)},
$$
where $N$ is the constant from \eqref{1.28.6}.
Similarly one takes care of the  term with $ds$
containing the derivatives of $\varphi_{m}$
in the second integral on the right in
\eqref{1.28.5}. Observing that by the dominated convergence
theorem and again by \eqref{1.27.9}
$$
\int_{0}^{T}|(\bar{\pi}_{s},|\beta^{k}_{s}|\,|\varphi_{m}-1|)\,ds
\to0\quad\text{(a.s.)},
$$
we conclude that the usual integrals on
the right-hand sides of
\eqref{1.28.5}
converge as $n\to\infty$
uniformly in $t\in[0,T]$ to the usual integral
the right-hand side of \eqref{1.28.1} almost surely.
 
To show the  convergence of the stochastic integrals
in \eqref{1.28.5} to the stochastic
integral in \eqref{1.28.1} uniform in probability it suffices
(and is necessary) to show that the quadratic
variation of the differences converges to zero in probability.
The said quadratic
variation is obviously less than a constant times
$$
\sum_{k}\int_{0}^{T}(\bar{\pi}_{s},\Lambda^{k}
(\varphi_{m}-1))^{2}\,ds,
$$
which tends to zero (a.s.) by the same reasons as above.
Thus, indeed  the right-hand sides of
\eqref{1.28.5}
converge as $n\to\infty$
uniformly in $t\in[0,T]$ in probability
to the right-hand side of \eqref{1.28.1}.
The left-hand sides converge for all $t\in[0,T]$
(a.s.) by the monotone convergence theorem.
This proves (i).

Assertion (ii) easily follows from the continuity
of $(\bar{\pi}_{t},1)$, the continuity of
$\bar{\pi}_{t}$ as an $L_{p}$-valued process, 
  and Scheff\'e's lemma. 

In (iii) that $P_{t}[\beta]$ is bounded 
follows from the boundedness
of $\beta$. The stated measurability properties of 
$P_{t}[\beta]$ are obtained by a standard measure-theoretic
argument form the fact that if $f(t,x,y)=\alpha(t)\beta(x)
\gamma(y)$, where $\alpha,\beta,\gamma$ are smooth functions 
with compact support, then
$$
\int_{\bR^{d}}f(t,x,y_{t})\bar{\pi}_{t}(x)\,dx
=\alpha(t)\gamma(y_{t})\int_{\bR^{d}}
\beta(x)\bar{\pi}_{t}(x)\,dx
$$
possesses the measurability properties in (iii) since
the last factor is a continuous (a.s.) $\cF^{y}_{t}$-adapted
process.

To prove \eqref{1.28.4} it suffices to use 
\eqref{10.13.3} which implies that for each
$t\in[0,T]$ and $y\in\bR^{d_{1}-d}$ 
$$
E(B(t,x_{t},y)\mid \cF^{y}_{t})=(\bar{\pi}_{t},1)^{-1}
\int_{\bR^{d}}B(t,x,y)\bar{\pi}_{t}(x)\,dx\quad
\text{(a.s.)}
$$
and then plug in here $y_{t}$ in place of $y$
in the argument of $B$, which
is possible because  $B(t,x,y)$
is Lipschitz in $y$ (even in $(x,y)$). This finishes
proving assertion (iii).

In (iv) the fact that $\check{w}_{t}$ is $\cF^{y}_{t}$-measurable
easily follows from an equivalent formula for
$\check{w}_{t}$:
$$
\check{w}_{t}=\int_{0}^{t}\Psi(s,y_{s})\,dy_{s}
-\int_{0}^{t}P_{s}[\beta]\,ds,
$$
where all terms on the right are $\cF^{y}_{t}$-measurable.
Furthermore, $\check{w}_{t}$ turns out to be
 an $\cF^{y}_{t}$-martingale on $[0,T]$. To check this,
take any $\cF^{y}_{t}$-stopping time $\tau\leq T$
and notice that $\tau$ is also an $\cF_{t}$-stopping time,
so that
$$
E\check{w}_{\tau}=E\int_{0}^{\tau}
 (\beta_{t}(x_{t})-P_{t}[\beta])\,dt.
$$
By using \eqref{1.28.4} and the fact that, by definition,
$\{t<\tau\}\in \cF^{y}_{t}$
we see that the right-hand side equals
$$
E\int_{0}^{T}
I_{t<\tau}(\beta_{t}(x_{t})-P_{t}[\beta])\,dt
=\int_{0}^{T}EI_{t<\tau} \beta_{t}(x_{t})\,dt
-\int_{0}^{T}EI_{t<\tau}P_{t}[\beta] \,dt
$$
$$
=\int_{0}^{T}EI_{t<\tau} \beta_{t}(x_{t})\,dt
-\int_{0}^{T}EI_{t<\tau}\big(E( \beta_{t}(x_{t})
\mid\cF^{y}_{t})\big) \,dt=0.
$$
 Thus, $E\check{w}_{\tau}=0$ for any $\cF^{y}_{t}$-stopping
time $\tau\leq T$ which combined with the $\cF^{y}_{t}$-adaptedness
of $\check{w}_{t}$ and its continuity in $t$ is well known
to be equivalent to saying that $\check{w}_{t}$ is an
$\cF^{y}_{t}$-martingale on $[0,T]$. Its quadratic variation
can be evaluated as the limit of sums of products
of increments and is, obviously, equal to the quadratic
variation of $\tilde{w}_{t}$, which, as we have seen in the proof
of Lemma \ref{lemma 10.21.4}, is a Wiener process. Therefore,
the quadratic variation of $\check{w}_{t}$ is that of a Wiener
process and by L\'evy's theorem $\check{w}_{t}$ is itself
a Wiener process  with respect to $\cF^{y}_{t}$.
This proves assertion (iv).

In (v) inequality \eqref{1.28.03} follows from
\eqref{1.28.3}, the fact that $\beta$ is bounded,
and the well-known properties of exponential martingales.
To prove \eqref{1.28.3} observe that 
\eqref{1.28.1} in terms of $P_{t}[\beta]$ 
and $\check{w}^{k}_{t}$ is rewritten as
$$
d(\bar{\pi}_{t},1)=(\bar{\pi}_{t},\beta^{k}_{t})
\beta^{k}_{t}(x_{t})\,dt+(\bar{\pi}_{t},\beta^{k}_{t})
\,d\tilde{w}^{k}_{t}
$$
$$
=(\bar{\pi}_{t},1)P_{t}[\beta^{k}]\beta^{k}_{t}(x_{t})\,dt
+(\bar{\pi}_{t},1)P_{t}[\beta^{k}]\,d\tilde{w}^{k}_{t}
$$
$$
=(\bar{\pi}_{t},1)|P_{t}[\beta]|^{2}\,dt
+(\bar{\pi}_{t},1)P_{t}[\beta^{k}]\,d\check{w}^{k}_{t}.
$$
Hence, $(\bar{\pi}_{t},1)$ satisfies the linear equation 
$$
d(\bar{\pi}_{t},1)=(\bar{\pi}_{t},1)|P_{t}[\beta]|^{2}\,dt
+(\bar{\pi}_{t},1)P_{t}[\beta^{k}]\,d\check{w}^{k}_{t},
$$
the unique solution of which with initial data $(\bar{\pi}_{0},1)
=(\pi_{0},1)=1$ is known to be given by \eqref{1.28.3}.
The theorem is proved.

\mysection{Appendix}
                                           \label{section 1.31.3}
The setting in this section
is somewhat different from that
of Section \ref{section 4.17.1}.
Let $(\Omega,\cF,P)$ be a complete probability space
with an increasing filtration $\{\cF_{t},t\geq0\}$
of complete with respect to $(\cF,P)$ $\sigma$-fields
$\cF_{t}\subset\cF$. Denote $\cP$ the predictable
$\sigma$-field in $\Omega\times(0,\infty)$
associated with $\{\cF_{t}\}$. Let
 $w^{k}_{t}$, $k=1,2,...$, be independent one-dimensional
Wiener processes with respect to $\{\cF_{t}\}$.  

We take a stopping time $\tau$ and for $t\leq\tau$
  we are considering the following
equation in $\bR^{d}$
\begin{equation}
                                          \label{11.13.1}
du_{t}=(L_{t}u_{t}-\lambda u_{t}
+D_{i}f^{i}_{t}+f^{0}_{t})\,dt
+(\Lambda^{k}_{t}u_{t}+g^{k}_{t})\,dw^{k}_{t},
\end{equation}
  where $u_{t}=u_{t}(x)=u_{t}(\omega,x)$ is an unknown function,
$$
L_{t}\psi(x)=D_{j}\big(a^{ij}_{t}(x)D_{i}\psi(x)
+a^{j}_{t}(x)\psi(x)
\big)+b^{i}_{t}(x)D_{i}\psi(x)+c_{t}(x)\psi(x),
$$
$$
\Lambda^{k}_{t}\psi(x)=\sigma^{ik}_{t}(x)D_{i}\psi(x)
+\nu^{k}_{t}(x)\psi(x),
$$
the summation  convention
with respect to $i,j=1,...,d$ and $k=1,2,...$ is enforced
and detailed assumptions on the coefficients and the 
free terms will be given later.

Fix a number
$$
p\geq2
$$
and denote $L_{p}=L_{p}(\bR^{d})$.
We use the same notation $L_{p}$ for vector- and matrix-valued
or else
$\ell_{2}$-valued functions such as
$g_{t}=(g^{k}_{t})$ in \eqref{11.13.1}. For instance,
if $u(x)=(u^{1}(x),u^{2}(x),...)$ is 
an $\ell_{2}$-valued measurable function on $\bR^{d}$, then
$$
\|u\|^{p}_{L_{p}}=\int_{\bR^{d}}|u(x)|_{\ell_{2}}^{p}
\,dx
=\int_{\bR^{d}}\big(
\sum_{k=1}^{\infty}|u^{k}(x)|^{2}\big)^{p/2}
\,dx.
$$

As above
$$
D_{i}=\frac{\partial}{\partial x^{i}},\quad i=1,...,d,
\quad\Delta=D^{2}_{1}+...+D^{2}_{d}.
$$
By $Du$ and $D^{2}u$ we mean the gradient and the matrix
of second order derivatives with respect
to $x$ of a function $u$ on $\bR^{d}$.

As above, for $\gamma\in\bR$
 by $H^{\gamma}_{p}=(1-\Delta)^{-\gamma/2}
L_{p}$ we denote the space
of Bessel potentials. 
Observe a slight change of notation. 
Since we will always be dealing with
$\bR^{d}$ we drop this symbol in the notation like
$H^{\gamma}_{p}(\bR^{d})$.
Most often in this appendix
we will use $H^{\gamma}_{p}$ for $\gamma=0,1$ 
and use \eqref{1.16.2}.

If $\tau$ is a stopping time, then
$$
\bH^{\gamma}_{p}(\tau):=L_{p}(\opar 0,\tau\cbrk,\cP,
H^{\gamma}_{p}),\quad\bL_{p}(\tau)=\bH^{0}_{p}(\tau). 
$$
We also need the space $\cH^{1}_{p}(\tau)$,
which is the space of functions $u_{t}
=u_{t}(\omega,\cdot)$ on $\{(\omega,t):
0\leq t\leq\tau,t<\infty\}$ with values
in the space of generalized functions on $\bR^{d}$
having the following properties:

(i) For any $T\in[0,\infty)$, we have $u
\in \bH^{1}_{p}(\tau\wedge T)$ 
and $u_{0}\in L_{p}(\Omega,\cF_{0},L_{p})$;

(ii) There exist   $f^{i}\in \bL_{p}(\tau)$,
$i=0,...,d$ and $g=(g^{1},g^{2},...)\in \bL_{p}(\tau)$
such that
 for any $\varphi\in C^{\infty}_{0}$ with probability 1
for all finite $t\leq\tau$
we have
\begin{equation}
                                                 \label{1.2.1}
(u_{t},\varphi)=(u_{0},\varphi)+
\int_{0}^{t}\big(-(f^{i}_{s},D_{i}\varphi)
+(f^{0}_{s},\varphi)\big)\,ds
+\sum_{k=1}^{\infty}\int_{0}^{t}
(g^{k}_{s},\varphi)\,dw^{k}_{s}.
\end{equation}

The reader can find 
in \cite{K99} a discussion of (i) and (ii),
in particular, the fact that the series in \eqref{1.2.1}
converges uniformly in probability on every finite
subinterval of $[0,\tau)$. On the other hand,
it is worth saying that
the above introduced space $\cH^{1}_{p}(\tau)$
are not quite the same as in \cite{K99}. There are 
three differences. One is that there is a restriction
on $u_{0}$ in \cite{K99}. However the most important
spaces are
$\cH^{1}_{p,0}(\tau)$ which are defined as the subsets of
$\cH^{1}_{p}(\tau)$ consisting of functions with
$u_{0}=0$. All other elements of $\cH^{1}_{p}(\tau)$
are obtained by adding to an element of $\cH^{1}_{p,0}(\tau)$
an  appropriate continuation for $t>0$ of the initial data.
Another issue is that in \cite{K99}
we have $f^{i}=0$, $i=1,...,d$, 
and $f^{0}\in\bH^{-1}_{p}(\tau)$. Actually, this 
difference is fictitious  because one knows that
any $f\in H^{-1}_{p}$ 

(a) has the form $D_{i}f^{i}+f^{0}$
with $f^{j}\in L_{p}$ and
$$
\|f\|_{H^{-1}_{p}}\leq N\sum_{j=0}^{d}\|f^{j}\|_{L_{p}},
$$
where $N$ is independent of $f,f^{j}$, and on the other hand,

(b) for any $f\in H^{-1}_{p}$  there exist
$f^{j}\in L_{p}$ such that $f=D_{i}f^{i}+f^{0}$ and
$$
\sum_{j=0}^{d}\|f^{j}\|_{L_{p}}\leq N\|f\|_{H^{-1}_{p}},
$$
where $N$ is independent of $f$.

The third difference is that instead of (i) we require
$D^{2}u\in\bH^{-1}_{p}(\tau)$ in \cite{K99}. However,
as it follows from Theorem 3.7 of \cite{K99}
and the boundedness of the operator $D:L_{p}\to H^{-1}_{p}$,
this difference disappears if $\tau$ is a bounded stopping time.

To summarize,
 the spaces $\cH^{1}_{p,0}(\tau)$ introduced above and
in \cite{K99} coincide if $\tau$ is bounded
and we choose a particular 
representation
of the deterministic part of the stochastic differential
just for convenience.  

In case that property (ii) holds, we write
\begin{equation}
                                       \label{12.3.1}
du_{t}=(D_{i}f^{i}_{t}+f^{0}_{t})\,dt
+g^{k}_{t}\,dw^{k}_{t}
\end{equation}
for $t\leq\tau$
and this explains the sense in which equation
\eqref{11.13.1} is understood. Of course, we still need to
specify appropriate assumptions on the coefficients
and the free terms in \eqref{11.13.1}.
Before we go to these assumptions we remind the reader that
according to \cite{K99} and the above discussion,
for bounded $\tau$,
one introduces a norm in
$\cH^{1}_{p,0}(\tau)$ by
$$
\|u\|_{\cH^{1}_{p,0}(\tau)}
=E\int_{0}^{\tau}\big(
\sum_{j=1}^{d}\|D_{j}u_{t}\|_{L_{p}}^{p}+
\sum_{j=0}^{d}\|f^{j}_{t}\|_{L_{p}}^{p}
+\|g_{t}\|_{L_{p}}^{p}\big)\,dt
$$
if $u$ satisfies \eqref{12.3.1}. By identifying two elements
of $\cH^{1}_{p,0}(\tau)$ if their difference has a zero 
$\cH^{1}_{p,0}(\tau)$-norm, 
one obtains a Banach space (see \cite{K99}).

We will also identify two elements 
$u',u'' \in\cH^{1}_{p}
(\tau)$ if and only if the difference $u'-u''$ is in
$\cH^{1}_{p,0}(\tau)$ and equals zero.

\begin{assumption}
                                        \label{assumption 1.2.1}

(i) The coefficients $a^{ij}_{t}$, $a^{i}_{t}$, $b^{i}_{t}$,
$\sigma^{ik}_{t}$, $c_{t}$, and $\nu^{k}_{t}$ are measurable 
with respect to $\cP\times B(\bR^{d})$, where $B(\bR^{d})$
is the Borel $\sigma$-field on $\bR^{d}$.

(ii) There is a constant $K$ such that
for all values of indices and arguments
$$
|a^{i}_{t}|+|b^{i}_{t}|+|c_{t}|+|\nu|_{\ell_{2}}\leq
K,\quad c_{t}\leq0.
$$

(iii) There is a constant $\delta>0$ such that
for all values of the arguments and $\xi\in\bR^{d}$
\begin{equation}
                                             \label{1.3.2}
(a^{ij}_{t}-\alpha^{ij}_{t})\xi^{i}
\xi^{j}\geq\delta|\xi|^{2},\quad|a^{ij}_{t}|\leq\delta^{-1},
\end{equation}
where   $\alpha^{ij}_{t}=(1/2)(\sigma^{i\cdot},\sigma^{j\cdot})
_{\ell_{2}}$. Finally, the constant $\lambda\geq0$.

\end{assumption} 

Assumption \ref{assumption 1.2.1} (i) guarantees
that equation \eqref{11.13.1} makes perfect sense
 for any constant $\lambda$
if $u\in\cH^{1}_{p}(\tau)$.
By the way, adding the term $-\lambda u_{t}$
with constant $\lambda\geq0$ is one more technically
convenient step. One can always introduce this term,
if originally it is absent, by considering $v_{t}:=u_{t}
e^{ \lambda t}$.

\begin{assumption}
                                       \label{assumption 1.2.6}
There is a continuous function $\kappa(\varepsilon)$
defined for $\varepsilon\geq0$ such that $\kappa(0)=0$
and 
$$
|\sigma^{i\cdot}_{t}(x)-
\sigma^{i\cdot}_{t}(x)|_{\ell_{2}}+
|a^{ij}_{t}(x)-a^{ij}_{t}(y)|\leq\kappa(|x-y|)
$$
for all $i,j,t,x,y$.
\end{assumption}

Here are the main results used in the previous sections
concerning \eqref{11.13.1}. They are taken from 
  \cite{Ki} and \cite{Kr09}. 
Generalization of these results to the case
of VMO coefficients $a^{ij}_{t}$ can be found in \cite{Kr09}. 

\begin{theorem}
                                    \label{theorem 12.7.1}
Let   $\lambda\geq0$, let
  $\tau$ be a stopping time, let
$f^{j},g\in\bL_{p}(\tau)$, and let
$\psi$ be a function such that $\psi\in\cH^{1}_{p}(\tau)
\cap\bH^{1}_{p}(\tau)$. Then
equation \eqref{11.13.1}
on $[0,\tau)$
has a unique solution  
$u\in\cH^{1}_{p}(\tau)$ such that
$u_{0}=\psi_{0}$.
 
Write
$$
d\psi_{t}=(D_{i}\alpha^{i}_{t}+\alpha^{0}_{t})\,dt
+\beta^{k}_{t}\,dw^{k}_{t}.
$$
Then the above  solution $u$   satisfies
$$
 \lambda^{1/2}\|u\|_{\bL_{p}(\tau)} 
+\|Du \|_{\bL_{p}(\tau)} 
$$
$$
  \leq N\big(\sum_{i=1}^{d}
\|f^{i}\|_{\bL_{p}(\tau)}+\|g\|_{\bL_{p}(\tau)}
+\sum_{i=1}^{d}
\|\alpha^{i}\|_{\bL_{p}(\tau)}+\|\beta\|_{\bL_{p}(\tau)}
+\| \psi\|_{\bH^{1}_{p}(\tau)}
\big)
$$
\begin{equation}
                                             \label{1.4.2}
+N\lambda^{-1/2}(\|f^{0}\|_{\bL_{p}(\tau)}
+\|\alpha^{0}\|_{\bL_{p}(\tau)}+\| \psi\|_{\bH^{1}_{p}(\tau)}
 )
+N\lambda^{1/2}\| \psi\|_{\bL_{p}(\tau)},
\end{equation}
provided that $\lambda>\lambda_{0}$, where
the constants $N,\lambda_{0}\geq0$
depend only on $d,p,K,\delta$, and the function
$\kappa$.
 
\end{theorem}

Observe that 
estimate \eqref{1.4.2} shows a good reason for
writing the free term in \eqref{11.13.1}
in the form $D_{i}f^{i}+f^{0}$, because
$f^{i}$, $i=1,...,d$, and $f^{0}$
enter \eqref{1.4.2} differently.

Here is a result about continuous dependence
of solutions on the data.

\begin{theorem}
                                    \label{theorem 1.4.1}
Assume that for each $n=1,2,...$
we are given functions $a^{nij}_{t}$, $a^{ni}_{t}$, $
b^{ni}_{t}$, $c^{n}_{t}$, $\sigma^{nik}_{t}$, $\nu^{k}_{t}$, 
$f^{ni}_{t}$, $g^{nk}_{t}$, and $\psi^{n}$
having the same meaning  and satisfying the same assumptions
with the same $\delta,K,\kappa$ as the original ones.
Assume that  
$$
(a^{nij}_{t},a^{ni}_{t},b^{ni}_{t},c^{n}_{t})\to
(a^{ ij}_{t},a^{ i}_{t},b^{ i}_{t},c _{t}),
$$
$$
|\sigma^{ni\cdot}_{t}-\sigma^{ i\cdot}_{t}|_{\ell_{2}}+
|\nu^{n}_{t}-\nu _{t}|_{\ell_{2}}\to0
$$
as $n\to\infty$ in measure $P(d\omega)dtdx$. Also let
$$
d\psi^{n}_{t}=(D_{i}\alpha^{ni}_{t}+\alpha^{n0}_{t})
\,dt+\beta^{nk}_{t}\,dw^{k}_{t}
$$
and assume that for a stopping time $\tau$
$$
\sum_{j=0}^{d}(\|f^{nj}-f^{j}\|_{\bL_{p}(\tau)}+
\|\alpha^{nj}-\alpha^{j}\|_{\bL_{p}(\tau)})
$$
$$
+
\|g^{n }-g\|_{\bL_{p}(\tau)}+\|\beta^{n }-\beta\|_{\bL_{p}(\tau)}
+\|\psi^{n}-\psi\|_{\bH^{1}_{p}(\tau)}\to0
$$
as $n\to\infty$. Take $\lambda\geq\lambda_{0}$, take
 the function $u$ from 
Theorem \ref{theorem 12.7.1} and let $u^{n}$ be unique solutions
of equations \eqref{11.13.1} constructed from
$a^{nij}_{t}$, $a^{ni}_{t}$, $
b^{ni}_{t}$, $c^{n}_{t}$, $\sigma^{nik}_{t}$, $\nu^{k}_{t}$, 
$f^{ni}_{t}$, and $g^{nk}_{t}$ and having initial
values  $\psi^{n}_{0}$. 

Then for any finite $T\geq0$ we have
$$
\|u^{n}-u\|_{\bH^{1}_{p}(\tau\wedge T)}\to0,\quad
E\sup_{t\leq\tau\wedge T}\|u^{n}_{t}-u_{t}\|_{L_{p}}^{p}\to0
$$ as $n
\to\infty$.

\end{theorem}

The following result shows that the solution
does not depend on $p$.

\begin{theorem}
                                    \label{theorem 1.4.2}

Let $p_{1},p_{2}\in[2,\infty)$ and let the assumptions
of Theorem \ref{theorem 12.7.1} be satisfied
with $p=p_{1}$ and $p=p_{2}$. Then the solutions corresponding
to $p=p_{1}$ and $p=p_{2}$ coincide, that is there
is a unique solution $u\in \cH^{1}_{p_{1}}(\tau)
\cap\cH^{1}_{p_{2}}(\tau)$
of equation \eqref{11.13.1} with initial data $\psi_{0}$.

\end{theorem}

In many situation the following maximum principle
is useful.

\begin{theorem}
                                    \label{theorem 1.4.3}

Under the assumptions
of Theorem \ref{theorem 12.7.1} suppose that $\psi_{0}\geq0$,
$f^{i}=0$, $i=1,...,d$, $f^{0}\geq0$, $g=0$.
Then for the solution $u$ almost surely
we have $u_{t}\geq0$ for all finite $t\leq\tau$.
 
\end{theorem}

Finally, we used the following embedding theorem
(see Corollary 4.12 and Remark 4.14 of \cite{Kr01}).
 For $\kappa\in(0,1)$, 
 a Banach space $X$,
and a set $A\subset\bR^{d}$
by $C^{\kappa}(A,X)$ we mean H\"older's space
of continuous $X$-valued functions on $A$ with finite
norm $\|\cdot\|_{C^{\kappa}(A,X)}$ defined by
$$
[|u|]_{C^{\kappa}(A,X)}=\sup_{s,t\in A}
|t-s|^{-\kappa}|u(t)-u(s)|_{X},\quad
\|u\|_{C (A,X)}=\sup_{ t\in A}
 |u(t) |_{X},
$$
$$
\|u\|_{C^{\kappa}(A,X)}=
[|u|]_{C^{\kappa}(A,X)}+\|u\|_{C (A,X)}.
$$

\begin{theorem}
                                         \label{theorem 1.17.1}
Let $\tau\leq T$, where the constant $T\in(0,\infty)$
and let $u\in\cH^{1}_{p}(\tau)$ satisfy \eqref{12.3.1}
with $f^{j}\in \bL_{p}(\tau)$, $g\in\bL_{p}(\tau)$,
and $u_{0}\in L_{p}(\Omega,\cF_{0}, H^{1-2/p}_{p})$,
Then:

(i) Almost surely   $u_{t}$ is a continuous
function of $t$ with values in
 $L_{p}$ for all  $t\in[0,\tau]$.

(ii) (case $p>2$)
Assume that for some numbers $\alpha$ and $\beta$ we have
$$
2/p<\alpha<\beta\leq1.
$$
Then, for any $a>0$,
\begin{equation}
                                                  \label{1.17.3}
E[u]^{p}_{C^{\alpha/2-1/p}([0,\tau],H^{1-\beta}_{p})}
\leq NT^{(\beta-\alpha)/p}a^{\beta-1}I(a),
\end{equation}
\begin{equation}
                                                  \label{1.17.4}
E\| u \|_{C 
([0,\tau],H^{1- \beta }_{p})}^{p}\leq 
N E\| u_{0}\|^{p}_{ H^{1-\beta}
_{p}}
+N  T^{p \beta /2-1} a^{ \beta -1}I(a),
\end{equation}
where the constants $N$ are independent of $a$,  
 $\tau$, $T$, and $u$ and
$$
I(a):=
a \|u\|^{p}_{\bH^{1}_{p }(\tau)}
+a^{-1}\|D_{i}f^{i}+f^{0}\|^{p}_{\bH^{-1}_{p }(\tau)}
+ \|g\|^{p}_{\bL_{p }(\tau )}.
$$
In particular, if $p(1-\beta)>d$, then 
\begin{equation}
                                                  \label{1.17.5}
E\sup_{x}[u(\cdot,x)
]^{p}_{C^{\alpha/2-1/p}([0,\tau])}
\leq NT^{(\beta-\alpha)/p}a^{\beta-1}I(a),
\end{equation}
\begin{equation}
                                                  \label{1.17.6}
E\sup_{t\in[0,T]}\| u(t,\cdot) \|_{C^{1-\beta-d/p}} ^{p}\leq 
N E\| u(0)\|^{p}_{ H^{1-\beta}
_{p}}
+N  T^{p \beta /2-1} a^{ \beta -1}I(a).
\end{equation}

Finally, \eqref{1.17.4} also holds if $p=2$ and $\beta=1$.

\end{theorem}

It is probably worth saying that \eqref{1.17.5} 
and \eqref{1.17.6} are not stated in \cite{Kr01}.
These are just obvious consequences
of \eqref{1.17.3} and \eqref{1.17.4} and the embedding theorem:
$H^{\gamma}_{p}
\subset C^{\gamma-d/p}$ if $\gamma-d/p>0$ and
$\gamma-d/p$ is not an integer.

\end{document}